\documentclass[11pt,a4paper]{article}
\usepackage{amsmath}
\usepackage{amsfonts}
\usepackage{amssymb}
\usepackage{amsthm}
\usepackage{mathrsfs}
\usepackage{dsfont}
\usepackage{paralist}
\usepackage{natbib}
\usepackage{bm}
\usepackage{color}
\usepackage[colorlinks=true,linkcolor=blue,citecolor=blue,pdfborder={0 0 0}]{hyperref}
\usepackage{graphicx}
\usepackage{tabularx}
\usepackage[left=3.5cm,right=3.5cm,top=2cm,bottom=2cm,includeheadfoot]{geometry}
\usepackage{comment}

\newtheoremstyle{normal}
{2ex}               
{3ex}               
{}                  
{}                  
{\bfseries} 
{}                  
{2pt}   
{\thmname{#1}\thmnumber{ #2.} \thmnote{(#3)}}

\newtheoremstyle{italic}
{2ex}
{3ex}
{\itshape}
{}
{\bfseries} 
{}
{2pt}
{\thmname{#1}\thmnumber{ #2.} \thmnote{(#3)}}

\theoremstyle{normal}
\newtheorem{definition}{Definition}[section]
\newtheorem{remark}[definition]{Remark}

\theoremstyle{italic}
\newtheorem{theorem}[definition]{Theorem}
\newtheorem{lemma}[definition]{Lemma}

\renewcommand{\P}{\mathbb{P}}
\newcommand{\F}{\mathcal{F}}

\newcommand{\E}{\mathbb{E}}
\newcommand{\N}{\mathds{N}}

\DeclareMathAccent{\verywidehat}{\mathord}{largesymbols}{'144}

\newcommand{\PP}{\mathbb{P}}

\newcommand{\Z}{\mathbb{Z}}

\newcommand{\Q}{\mathbb{Q}}

\newcommand{\Ga}{\Gamma}
\newcommand{\si}{\sigma}

\newcommand{\Om}{\Omega}

\newcommand{\f}{\mathcal F}

\renewcommand{\lll}{\mathcal L}

\newcommand{\proba}{(\Omega ,\f,(\f_t)_{t\geq0},\P)}

\newcommand{\probp}{(\Omega' ,\f',(\f'_t)_{t\geq0},\PP')}
\def\probt{\hbox{$(\widetilde{\Omega},\widetilde{\f},(\widetilde{\f}_t)_{t\geq0},\widetilde{\PP})$}}

\newcommand{\tols}{~\stackrel{\lll-(s)}{\longrightarrow}~}
\newcommand{\tol}{\stackrel{\lll}{\longrightarrow}}

\newcommand{\pn}{\stackrel{\P}{\longrightarrow}}

\allowdisplaybreaks[1]

\begin{document}

\title{A note on central limit theorems for quadratic variation in case of endogenous observation times
}

%
%
%
%

\author{Mathias Vetter\thanks{Christian-Albrechts-Universit\"at zu Kiel, Mathematisches Seminar, Ludewig-Meyn-Str.\ 4, 24118 Kiel, Germany.
{E-mail:} vetter@math.uni-kiel.de} \and Tobias Zwingmann\thanks{Philipps-Universit\"at Marburg, Fachbereich Mathematik und Informatik, Hans-Meerwein-Str., 35043 Marburg, Germany.
{E-mail:} zwingmann@mathematik.uni-marburg.de}}

\maketitle

\begin{abstract}
This paper is concerned with a central limit theorem for quadratic variation when observations come as exit times from a regular grid. We discuss the special case of a semimartingale with deterministic characteristics and finite activity jumps in detail and illustrate technical issues in more general situations. 
\end{abstract}



\smallskip



\section{Introduction}
\def\theequation{1.\arabic{equation}}
\setcounter{equation}{0}

High-frequency statistics has attracted a lot of attention in recent years. Given observations of a semimartingale $X$, one is often interested in estimation of its quadratic variation (or parts thereof, such as integrated volatility) and with associated central limit theorems. 

A natural way is to work in a setting where observations come at regular times, that is we have data $X_{j/n}$, $j=0, \ldots, n$, over the interval $[0,1]$, say. The most general paper on asymptotics in this setting is \cite{jacod2008} where various (stable) central limit theorems for functionals of discretely observed It\=o semimartingales are stated, including those for realized variance $RV(X,X)_t^n$ with 
\[
RV(X,Y)_t^n = \sum_{j=1}^{\lfloor nt \rfloor} (X_{j/n} - X_{(j-1)/n})(Y_{j/n} - Y_{(j-1)/n})
\]
for arbitrary processes $X$ and $Y$; see also \cite{jacpro1998} for earlier results on related statistics in the case of L\'evy processes. Suppose, $X$ is defined on the filtered probability space $\proba$ and can be decomposed as
\begin{align}\label{sm2}
X_t = X_0 + \int_0^t a_s ds + \int_0^t \sigma_s dW_s&+\int_0^t\int_{\mathds{R}}\kappa(\delta(s,x))(\mu-\nu)(ds,dx)\\
&\notag+\int_0^t\int_{\mathds{R}}\bar\kappa(\delta(s,x))\mu(ds,dx),
\end{align}
where $a$, $\sigma$ and $\delta$ are adapted processes and a truncation function $\kappa$, $\bar\kappa(x)=x-\kappa(x)$, separates large from compensated small jumps. The compensating intensity measure $\nu$ of the Poisson random measure $\mu$ admits the form $\nu(ds,dx)=ds\otimes \lambda(dx)$ for a $\sigma$-finite measure $\lambda$. In this case, the quadratic variation process becomes 
\[
[X,X]_t = \int_0^t \sigma_s dW_s + \sum_{0 \leq s \leq t} \Delta X^2_s, \qquad \Delta X_s = X_s - X_{s-}
\]
 Under essentially no extra conditions, Theorem 2.11 in \cite{jacod2008} gives the $\F$-stable central limit theorem 
\begin{align}
Z_t^n(X) = \sqrt{n} \left(RV(X,X)_t^n - [X,X]_t \right) \tols Z_t, \label{thmjac}
\end{align}
pointwise in $t$, where the limiting process $Z_t = U_t + V_t$ consists of two parts: The first one,
\[
U_t = \sqrt 2 \int_0^t \sigma_s^2 dW'_s,
\]
is due to the continuous martingale part of $X$ only. Here, $W'$ denotes an independent Brownian motion defined on an extension $\probp$ of the original space. The process $V_t$ takes a more complicated form and comes from both the continuous and the jump part of $X$. It is given by
\[
V_t = 2 \sum_{S_p \leq t} \Delta X_{S_p} \left(\kappa_p \sigma_{S_p-} R_p + \sqrt{1- \kappa_p} \sigma_{S_p} R'_p \right),
\]
where the sequence $(S_p)_p$ denotes an enumeration of the jump times of $X$, and where the $R_p$ and $R'_p$ are standard normal and the $\kappa_p$ are uniform random variables on $[0,1]$, all defined on the same extension as $W'$ and all mutually independent. In the case of a continuous $\si$, the limit obviously reduces to 
\begin{align} \label{limitV}
V_t = 2 \sum_{S_p \leq t} \Delta X_{S_p} \sigma_{S_p} R_p.
\end{align}
For details on stable convergence we refer to Section VIII.5 of \cite{jacshi2003}. 

Over the last decade, a lot of work has been connected with extensions of \eqref{thmjac} whenever the ideal setting of observations of $X$ at equidistant times is not realistic. Typically, the focus has been on a version of $X$ with continuous paths, which means that \eqref{sm2} reduces to 
\begin{align} \label{sm1}
X_t = X_0 + \int_0^t a_s ds + \int_0^t \sigma_s dW_s.
\end{align}
Regarding non-equidistant observations, authors have typically worked in settings where observations times come either in a deterministic way or are (essentially) independent of $X$. In this case, one still has stable convergence of the standardized realized variance to a limiting process $U'_t$, but it differs from $U_t$ by an extra factor which accounts for the variability of time. See for example \cite{hayetal2011}, \cite{mykzha2012} or \cite{koike2014}. A similar result can be obtained if jumps are present, but in this case with a different process $V'_t$ as well. See \cite{bibvet2015}. 

Even though important from a practical point a view, the situation with endogenous observation times has found much less attention. This is certainly the case because the proofs of the corresponding results become much more complicated then. For model (\ref{sm1}), central limit theorems are provided in \cite{fukasawa2010} and \cite{lietal2014}, once certain extra conditions are satisfied which might be difficult to check in practice. As a specific example, hitting times of a grid are discussed in \cite{fukros2012}. 

However, to the best of our knowledge no results exist in the general model (\ref{sm2}) involving jumps. The aim of this work therefore is to shed some light on this issue. In particular, we will prove a result similar to (\ref{thmjac}), but only in the less general case of a continuous It\=o semimartingales plus finite activity jumps with deterministic characteristics, observed at hitting times of a regular grid. This is a bit unsatisfactory from a practical point of view, but we will discuss the reasons for this slight limitation and possible guidelines for future research in order to solve the problem in a general framework. 

This work is organised as follows: Section \ref{stand} gives a short review on the case of deterministic observation times, and we give a heuristic explanation why the limiting process $Z_t$ is of the form as stated in (\ref{thmjac}). Section \ref{endog} contains the main theorem of this work, as well as a comparison with the standard result and remarks on further extensions of the model. All proofs are gathered in Section \ref{app}. 

\section{The standard case: Deterministic observations} \label{stand}
\def\theequation{2.\arabic{equation}}
\setcounter{equation}{0}

Let us shortly sketch the strategy which leads to the main result (\ref{thmjac}) from \cite{jacod2008}. The key idea is to decompose $X$ for each integer $q$ as 
\begin{align} \label{decomp}
X_t = C_t + N(q)_t + X(q)_t,
\end{align}
where
\begin{align*}
N(q)_t = \int_0^t \int \delta(s,x) 1_{\{\gamma(x) > 1/q\}} \mu(ds,dx)
\end{align*}
denotes the large jumps of $X$,
\[
C_t = \int_0^t \sigma_s dW_s
\] 
corresponds to the continuous martingale part of $X$, and the remainder $X(q)_t$ involving drift and small jumps is defined implicitly. One uses here the mild assumption that the process $|\delta(s,x)|$, which is responsible for the jump sizes, is bounded by some deterministic function $\gamma(x)$. 

Using (\ref{decomp}) and the binomial theorem one obtains
\begin{align*}
Z_t^n(X)
&= \sqrt{n} \left(RV(C,C)_t^n - [C,C]_t \right) + \sqrt{n} \left(RV(N(q),N(q))_t^n - [N(q),N(q)]_t \right) \\ &+ \sqrt{n} \left(RV(X(q),X(q))_t^n - [X(q),X(q)]_t \right) \\ &+ 2 \sqrt{n} RV(C,N(q))_t^n + 2 \sqrt{n}RV(C,X(q))_t^n + 2 \sqrt{n} RV(N(q), X(q))_t^n 
\end{align*}
for each integer $q$. A huge part of the proof of (\ref{thmjac}) deals with negligibility of most terms in the decomposition above. This is the case for all terms involving the remainder $X(q)$ which can be shown to converge to zero in probability if we first let $n \to \infty$ and then $q \to \infty$. The proof relies heavily on martingale techniques when dealing with the compensated small jumps of $X$. See for example Appendix B in \cite{bibvet2015}. Moreover, since $N(q)$ is a finite activity jump process, for any fixed $q$ we have 
\[
\sqrt{n} \left(RV(N(q),N(q))_t^n - [N(q),N(q)]_t \right) = 0
\]
identically, with a probability converging to one, as well. This explains why only two processes appear in the limit in (\ref{thmjac}). In case of continuous paths, the relevant term is $ \sqrt{n} \left(RV(C,C)_t^n - [C,C]_t \right)$, for which the proof of stable convergence to $U_t$ is very well understood by now. The main condition, due to conditional Gaussianity, is the convergence of the empirical conditional variance
\begin{align} \label{convcond} \nonumber
& \frac{2n}{3} \sum_{j=1}^{\lfloor nt \rfloor} \mathbb E\left[ (C_{j/n} - C_{(j-1)/n})^4 |\mathcal F_{(j-1)/n} \right] \\ =& \frac{2n}{3} \sum_{j=1}^{\lfloor nt \rfloor} \si^4_{(j-1)/n}\mathbb E\left[ (W_{j/n} - W_{(j-1)/n})^4\right] + o_{\P}(1)\stackrel{\mathbb P}{\longrightarrow} 2 \int_0^t \sigma_s^4 ds,
\end{align}
which determines the distribution of the limiting Brownian martingale. 

Let us therefore focus on the mixed part. Since $N(q)_t$ is a finite activity jump process, we can write 
\begin{align*}
2 \sqrt{n} RV(C,N(q))_t^n = 2 \sqrt n \sum_{S_p \leq t} \Delta N(q)_{S_p} \left(C_{i_+(S_p)/n} - C_{i_-(S_p)/n}\right),
\end{align*}
where $i_+(s)$ and $i_-(s)$ denote the index of the first observation past or equal to $s$ and the last observation prior to $s$, respectively. Even though the process looks already like $V_t$, there are several steps necessary in order to get to the final result. First, one uses a similar discretization argument as for the conditional variance in order to write 
\[
\left(C_{i_+(S_p)/n} - C_{i_-(S_p)/n}\right) = \sigma_{i_-(S_p)/n} \left(W_{S_p} - W_{i_-(S_p)/n}\right) + \sigma_{S_p} \left(W_{i_+(S_p)/n} - W_{S_p}\right),
\]
up to an error of size $o_{\P}(n^{-1/2}).$ Second, for each $p$ we have the equality 
\begin{align}
\label{step1} \left(W_{S_p} - W_{\frac{i_-(S_p)}n}, W_{\frac{i_+(S_p)}n} - W_{S_p}\right) = \left(\sqrt{S_p - \frac{i_-(S_p)}n} R_p, \sqrt{\frac{i_+(S_p)}n - S_p} R'_p \right)
\end{align}
in distribution. Third, the stable convergence 
\begin{align} \label{step2}
n\left(S_p - (i_-(S_p))/n\right) \tols \kappa_p 
\end{align}
follows, since the jump times of $N(q)$ are uniformly distributed over [0,1]; see the proof of Lemma 6.2 in \cite{jacpro1998}. All of these steps can be shown to hold jointly for $1 \leq p \leq k$, for an arbitrary integer $k$. This is sufficient to prove 
\[
2 \sqrt{n} RV(C,N(q))_t^n \tols 2 \sum_{S_p \leq t} \Delta N(q)_{S_p} \left(\kappa_p \sigma_{S_p-} R_p + \sqrt{1- \kappa_p} \sigma_{S_p} R'_p \right)
\]
for any fixed $q$. Letting $q \to \infty$, stable convergence to $V_t$ follows. Finally, some more technicalities are needed in order to prove the joint stable convergence of the two sequences leading to $U_t$ and $V_t$.

\section{The endogenous case: Exit times from a regular grid} \label{endog}
\def\theequation{3.\arabic{equation}}
\setcounter{equation}{0}

What makes an extension of (\ref{thmjac}) to random observation times $\tau_j^n$ difficult in general, is in particular the equality in distribution stated in (\ref{step1}). As long as the times are exogenous, i.e.\ independent of $X$, or satisfy some predictability condition, this equality still holds because at least the local behaviour of $W$ is independent of the observation times. The central limit theorem for 
\[
RV(X,X)_t^n = \sum_{0 < \tau_j^n \leq t} (X_{\tau_j^n} - X_{\tau_{j-1}^n})^2, \qquad \tau_0^n = 0, 
\]
then looks similar to (\ref{thmjac}), with some minor modifications because one has to find suitable conditions under which a result like (\ref{step2}) holds. See e.g.\ \cite{bibvet2015} for details. 

\begin{figure}
	\centering
	\includegraphics[width=\linewidth]{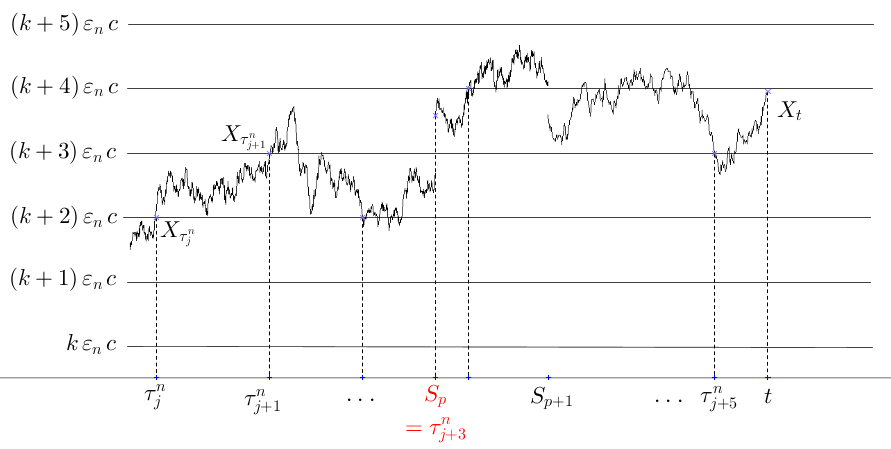}
	\caption[The endogenous observation scheme]{\small The observation scheme used: The first observations $\tau_j^n$ to $\tau_{j+2}^n$ are due to movements of the underlying Brownian motion only; the jump at $S_p$ is big enough to cause the observation $\tau_{j+3}^n$. Remark that between $\tau_{j+2}^n$ and $\tau_{j+3}^n$ the path hits $(k+2)\epsilon_n c$ several times without generating an observation. Last, the jump at $S_{p+1}$ is too small for the path to exit the interval $[(k+5)\epsilon_n c, (k+3)\epsilon_n c]$.}
	\label{fig:model}
\end{figure}

In the purely endogenous case, however, the observation times depend strongly on the process $X$ (and in particular on the process $W$), and therefore one cannot expect a version of (\ref{step1}) to remain valid. For this reason, we will work in a specific setting, which makes a computation of the distribution of increments of the Brownian part between successive observations possible. Suppose therefore in the following that we observe
\begin{align} \label{sm3}
X_t = \int_0^t a_s ds + \int_0^t \sigma_s dW_s + J_t,
\end{align}
where $a_s$ and $\sigma_s$ are deterministic continuous processes, $W_t$ is a standard Brownian motion and 
$
J_t
$
is a finite activity jump process with a deterministic compensator. We also assume that the function $\si$ is positive everywhere. Observations are coming as exit times from a regular grid. That is, for a given constant $c > 0 $ and a sequence $\epsilon_n \to 0$, which governs the asymptotics, we observe $X$ at the stopping times $\tau_0^n=0$ and 
\[\tau_j^n = \inf[ t > \tau_{j-1}^n: X_t \notin A_{X_{\tau_{j-1}^n}} ],\] 
where either \[ A_y = [y-c \epsilon_n, y + c \epsilon_n]\] if $y = k c \epsilon_n$ for some $k \in \mathbb Z$ or \[A_y = [\lfloor (c \epsilon_n)^{-1}y \rfloor c \epsilon_n, (\lfloor (c \epsilon_n)^{-1}y \rfloor +1) c \epsilon_n]\] if no such $k$ exists.  
A sketch of this observation scheme can be found in Figure~\ref{fig:model}. 

Roughly speaking, data is collected whenever the unit of the price is changing, e.g.\ from one cent to another. No stopping, however, takes place at the previously reached barrier, which is necessary because Brownian motion hits it infinitely often after starting from there. This scheme is a particular case of the setting in \cite{fukros2012}, who derived a central limit theorem for realized variance of a \textit{continuous} semimartingale. Note that in their setting observations came at hitting times (not necessarily of a regular grid), which is essentially similar to exit times for processes with continuous paths. 

%
%

In order to state the stable central limit theorem, we again assume the process $X$ to live on the filtered probability space $\proba$, and we consider a second probability space $\probp$ which supports a Brownian motion $W'$ and a sequence $(\eta_p)_{p \geq 1}$ of random variables with density 
\begin{align} \label{dens}
h(y) = \int_0^\infty \frac 1{\sqrt{2 \pi z}} \sum_{m=-\infty}^\infty (e^{-\frac{(y-4m)^2}{2z}} - e^{-\frac{(y+2 + 4m)^2}{2 z}}) dz 1_{[-1,1]}(y).
\end{align}
All of these are mutually independent and independent of $W'$ as well. Limiting variables are then defined on the product $\probt$ of the two afore-mentioned probability spaces. 
 
Finally, we denote with $S_1, S_2, \ldots$ a specific sequence of stopping times exhausting the jumps of $X$, namely we order the sequence in chronological order. This is possible, as there are only finitely many jumps almost surely. In principle, other enumerations were possible as well, but this choice facilitates some parts of the proof.
 
\begin{theorem} \label{thm1}
Suppose that $X$ is given by \eqref{sm3} and is observed at the random instants coming from the sampling scheme described above. Then, as $n \to \infty$, we have the $\mathcal F$-stable convergence
\[
\widetilde Z_t^n(X) = \epsilon_n^{-1} \left( RV(X,X)_t^n - [X,X]_t \right) \tols \widetilde Z_t = \widetilde U_t + \widetilde V_t,
\] 
pointwise in $t$, where 
\[
\widetilde U_t = \sqrt{\frac 23} c \int_0^t \sigma_s dW'_s
\]
and 
\[
\widetilde V_t = 2 \sum_{S_p \leq t} \Delta X_{S_p} c \eta_p,
\]
with the Brownian motion $W'$ and the sequence $(\eta_p)_{p \geq 1}$ as above.
\end{theorem}

\begin{remark} 
\begin{itemize}
	\item[(i)] In contrast to (\ref{limitV}), the distribution of $\widetilde{V}_t$ does not involve a factor $\sigma_{S_p}$. Intuitively, when finding analogues of (\ref{step1}) and (\ref{step2}) we are concerned with the distribution of the random variables 
	\begin{align*}
	\epsilon_n^{-1} \si_{\tau_n^{-}(S_p)} \left(W_{S_p} - W_{\tau_n^{-}(S_p)}\right),
	\end{align*} 
	where it is sufficient to take $S_p$ as the right end point because a jump of $X$ causes an exit from the current interval with probability converging to one. (Here and below, $\tau_n^{-}(t)$ is the last observation prior to $t$.) After rescaling, each distribution essentially equals the one of a Brownian motion $\si W$ (starting in zero and with $\si = \si_{\tau_n^{-}(S_p)}$) at time $t$, when it is known to not having left the interval $[-c, c]$ by that time. This additional knowledge makes the distribution independent of $\si$. See Lemma \ref{step3} for details. 
	\item[(ii)] A similar effect is known from \cite{fukros2012}, from which we borrow the stable convergence of the continuous martingale part. Whereas in the case of exogeneous observation times $\sigma$ enters to second order in the conditional standard deviation, as seen e.g.\ from (\ref{convcond}), it only contributes to first order to the (conditional) standard deviation of $\widetilde U$. The missing order is due to a factor $c$ in our situation, which enters likewise in the mixed term.   
	\item[(iii)] We believe that a generalization to the case of random characteristics and irregular grids holds as well, but in this situation the proof becomes even more involved. If the grid takes the form
	\[
	\mathbb G_n = \{\epsilon_n p_i|i\in \mathbb Z\},
	\]
	we conjecture that apart from $\widetilde U$, which changes according to the results in \cite{fukros2012} under their assumptions, the second summand becomes
	\[
\widetilde V_t = 2 \sum_{S_p \leq t} \Delta X_{S_p} \eta_p(X_{S_p-})
	\] 
	where the distribution of $\eta_p(x)$ is a suitably weighted mixture of the limiting distribution of two independent Brownian motions starting at the grid points right above and right below $x$, respectively, and which are known to not having caused another observation. 
	\item[(iv)] Using our strategy of proof it is extremely difficult to prove (a version of) Lemma \ref{step3} and, thus, Theorem \ref{thm1} when one works with infinitely many jumps within $X$. In our setting, as two jumps are usually far apart, the key steps in the proof regard the exact distribution of exit times of Brownian motion and a result by \cite{freedman1983} who derives an explicit formula for 
\begin{align} \label{fried}
\P\left( W_z \leq  x, \inf_{0 \leq u \leq z} W_u > -c, \sup_{0 \leq u \leq z} W_u < -c \right).
\end{align}
In the general case, one has to account for exit times of jump processes as well, and very little is known even for L\'evy processes. (Asymptotics are provided in \cite{rostan2011}, however.) Also, to the best of the authors' knowledge, a generalization of (\ref{fried}) does not even exist in the literature for Brownian motions with a constant drift, let alone for more general processes involving jumps. \qed
\end{itemize}
\end{remark}

\section{Appendix: Proof of Theorem \ref{thm1}} \label{app}
\def\theequation{4.\arabic{equation}}
\setcounter{equation}{0}
We start with the proof of Theorem \ref{thm1} when $X_t = \si W_t + J_t$ with a constant $\si > 0$. As above, we set $C_t = \sigma W_t$, and we denote with $S_1, S_2, \ldots$ the sequence of jump times in chronological order. For convenience, $S_0=0$. We write
\[
U_t^n = \epsilon_n^{-1} \left( RV(C,C)_t^n - [C,C]_{\tau_n^{-}(t)} \right).
\]
Furthermore, we set $$\alpha(n,p) = \sigma \epsilon_n^{-1} \left( W_{S_p} - W_{\tau_n^{-}(S_p)} \right)$$ for each $p$. The main part of the proof is Theorem \ref{thm1} is contained in the following lemma. 

\begin{lemma} \label{lemstab}
As $n \to \infty$, we have the $\mathcal F$-stable convergence
\[
(U_t^n,(\alpha(n,p))_{p \geq 1}) \tols (\overline U_t,(\eta_p)_{p \geq 1}),
\]
with $\overline U_t = \sqrt{2/3} c \sigma W'_t$ and the sequence $\eta_p$ as before.
\end{lemma}

\textbf{Proof.} \rm Using the standard metric on an infinite Cartesian product, it is enough to prove the result for $1 \le p \le k$ and some $k$. Without loss of generality, we assume $k \leq \Gamma$ where $\Gamma$ is the random, almost surely finite number of jumps over $[0,t]$. 
Thus, we show 
\begin{align*} 
(U_t^n,\alpha(n,p)_{1 \leq p \leq k})) \tols (\overline U_t,(\eta_p)_{1 \leq p \leq k}), 
\end{align*}
which formally means that we have to prove 
\begin{align*}
\E\big[\Psi g(U_t^n) \prod_{p=1}^k h_p(\alpha(n,p))\big] \to \widetilde \E\big[\Psi g(\overline U_t) \prod_{p=1}^k h_p(\eta_p)] 
\end{align*}
for any bounded $\Psi \in \mathcal F$ and for all bounded Lipschitz functions $g,h_1, \ldots, h_k$. 


First, the idea is to separate $U_t^n$ from the other variables. To this end, we need some technical definitions. Fix some integer $\ell$ and set $S_p^{\ell-}=(S_p-1/\ell)_+$ and $S_p^{\ell+}=S_p+1/\ell$. Then $B_\ell = \cup_{p = 1}^\Gamma[S_p^{\ell-},S_p^{\ell+}]$ denotes the union of all those intervals around the jump times, and we use $\Lambda_n(\ell,t)$ to define the set of indices $j$ such that $\tau_j^n \leq t$ and $[\tau_{j-1}^n,\tau_j^n] \cap B_\ell = \emptyset$. The increments over all such $j$ then cover 
\[
D_n(\ell,t) = [0,t] \backslash \left( \cup_{p = 1}^\Gamma [\tau_n^-(S_p^{\ell-}),\tau_n^+(S_p^{\ell+})] \right).
\]
Finally, we define 
\[
U(\ell)_t^n = \epsilon_n^{-1} \sum_{j \in \Lambda_n(\ell,t)} \left( (C_{\tau_j^n} - C_{\tau_{j-1}^n})^2 - \sigma^2 (\tau_j^n - \tau_{j-1}^n) \right),
\]
which is the same quantity as $U_t^n$, but only involves intervals which are far from the jump times. We will see that $U(\ell)_t^n$ and $U_t^n$ obey the same asymptotic law.

It is helpful to introduce the sets 
\[
\Omega_n = \{ \omega \in \Omega: \inf\{ |\Delta J_{S_p}|; 1 \leq p \leq \Gamma \} > 2 \epsilon_n c \}
\]
and 
\begin{align*}
\Omega(\ell) = \left \{ \omega \in \Omega: \inf\{|S_p - S_{p-1}|; 1 \leq p \leq \Gamma \} > \frac 2 \ell \right \}.
\end{align*}
On the set $\Omega_n$ every jump of $X$ leads to an exit from the current interval, thus causes an observation.  Obviously, $\Omega_n \to \Omega$ as $n \to \infty$, because there are only finitely many jumps almost surely. Therefore, it is equally well possible to prove 
\begin{align*}
\E\big[\Psi g(U_t^n) \prod_{p=1}^k h_p(\alpha(n,p)) 1_{\Omega_n} \big] \to \widetilde \E\big[\Psi g(\overline U_t) \prod_{p=1}^k h_p(\eta_p)]. 
\end{align*}
Let us further remark that on $\Omega(\ell)$ any two jumps are further than $2/\ell$ apart, and thus $B_{\ell}$ becomes a disjoint union. 

We are now interested in showing
\begin{align} \label{steplem}
\sup_{n \in \N} \big|\E\big[\Psi (g(U_t^n) -  g(U(\ell)_t^n) \prod_{p=1}^k h_p(\alpha(n,p)) 1_{\Omega_n} \big] \big| \to 0
\end{align}
as $\ell \to \infty$, and because of $\Omega(\ell)  \to \Omega$ we may assume to live on $\Omega(\ell)$ as well. Let $(\mathcal F''_t)_t$ denote the filtration which is the smallest one containing $(\mathcal F_t)_t$ and such that the jump measure $\mu$ of $J$ is $\mathcal F''_0$-measurable. Note that $W$ remains a standard Brownian motion with respect to this filtration. We discuss
\begin{align*} 
\E\left[|U_t^n  - U(\ell)_t^n|^2 1_{\Omega_n} 1_{\Omega(\ell)} \right] = \E\left[1_{\Omega_n} 1_{\Omega(\ell)} \E \left[ |U_t^n  - U(\ell)_t^n|^2 \vert \mathcal F''_0 \right] \right]
\end{align*}
first. On the set $\Omega_n \cap \Omega(\ell)$ all jump times are observation times, and conditionally on $\mathcal F''_0$ these are known. Therefore, by successive conditioning and using the Burkholder-Davis-Gundy inequality,
\begin{align*}
& 1_{\Omega_n \cap \Omega(\ell)} \E \left[ |U_t^n  - U(\ell)_t^n|^2 \vert \mathcal F''_0 \right] \\
 =& 1_{\Omega_n \cap \Omega(\ell)} \epsilon_n^{-2} \E \big[ \big \vert  \sum_{j \notin \Lambda_n(\ell,t)} \big( (C_{\tau_j^n} - C_{\tau_{j-1}^n})^2 - \sigma^2 (\tau_j^n - \tau_{j-1}^n) \big) \big \vert^2  \vert \mathcal F''_0 \big] \\ =&1_{\Omega_n \cap \Omega(\ell)} \epsilon_n^{-2} \E \big[ \sum_{j \notin \Lambda_n(\ell,t)} \big( (C_{\tau_j^n} - C_{\tau_{j-1}^n})^2 - \sigma^2 (\tau_j^n - \tau_{j-1}^n) \big)^2  \vert \mathcal F''_0 \big] \\ \leq& K 1_{\Omega_n \cap \Omega(\ell)}  \epsilon_n^{-2}\E \big[ \sum_{j \notin \Lambda_n(\ell,t)} (C_{\tau_j^n} - C_{\tau_{j-1}^n})^4  \vert \mathcal F''_0 \big] \\ \leq& K 1_{\Omega_n \cap \Omega(\ell)}  \epsilon_n^{-2} \sum_{p=1}^\Gamma \E \big[ \sum_{j: \tau_n^-(S_p^{\ell-}) < \tau_j^n \leq \tau_n^+(S_p^{\ell+})} (C_{\tau_j^n} - C_{\tau_{j-1}^n})^4  \vert \mathcal F''_0 \big].
\end{align*}
Here and below, $K$ denotes an unspecified constant. Note that on $\Omega(\ell)$ all observations after $S_p^{\ell-}$ and until $S_p^{\ell+}$ are due to exits of the Brownian motion from the respective interval. The final observation $\tau_n^+(S_p^{\ell+})$, for example, might be due to the next jump $S_{p+1}$, however, but apart from these two additional increments -- which, as can be seen from (\ref{moeper}) and (\ref{renewal}) later, do not alter the asymptotic behaviour -- we are in the same situation as in \cite{fukros2012}. The proof of their Proposition 4.3 then gives
\[
\epsilon_n^{-2} \sum_{p=1}^\Gamma \E \big[ \sum_{j: \tau_n^-(S_p^{\ell-}) < \tau_j^n \leq \tau_n^+(S_p^{\ell+})} (C_{\tau_j^n} - C_{\tau_{j-1}^n})^4  \vert \mathcal F''_0 \big] \pn \si^2 c^2 \frac{2 \Ga}{\ell}, 
\] 
which basically equals the variance of realized variance when computed over intervals outside of $\Lambda_n(\ell,t)$ only. Then, it it easy to deduce
\[
\sup_{n \in \N} \big|\E\big[\Psi (g(U_t^n) -  g(U(\ell)_t^n) \prod_{p=1}^k h_p(\alpha(n,p)) 1_{\Omega_n} 1_{\Omega(\ell)} \big] \big| \leq K \Big(\frac{\E[\Ga]}{\ell}\Big)^{1/2},
\]
because the random variables inside the expectation are uniformly bounded by a constant, thus in $L^2$, and because it suffices to bound $|U_t^n  - U(\ell)_t^n|$, due to the Lipschitz property of $g$. Since $\E[\Gamma] < \infty$, we obtain (\ref{steplem}) as requested. Also, since $B_\ell$ decreases to a discrete set we have
\begin{align} \label{ell2}
\overline U(\ell)_t := \sqrt{\frac 23} c  \sigma  \int_0^t 1_{B_\ell^c}(s) dW'_s \pn \overline U_t,
\end{align}
as $\ell \to \infty$. 

Combining (\ref{steplem}) and (\ref{ell2}) it is certainly enough to show 
\begin{align} \label{step7}
\E\big[1_{\Omega_n \cap \Omega(\ell)} \Psi g(U(\ell)_t^n) \prod_{p=1}^k h_p(\alpha(n,p))\big] \to \widetilde \E\big[1_{\Omega(\ell)} \Psi g(\overline U(\ell)_t) \prod_{p=1}^k h_p(\eta_p)] 
\end{align} 
%
%
%
for any fixed integer $\ell$. Let us introduce 
\[
\Omega_n(\ell) = \bigcap_{1 \leq p \leq \Gamma} \left( \{\tau_n^+(S_p) < S_p^{\ell+} \} \cap \{\tau_n^-(S_p) > S_p^{\ell-} \} \right) \cap \Omega_n.
\] 
Obviously, 
\begin{align} \label{abschP}
\P \big(\Omega_n(\ell)^c \big) \leq \E\big[\sum_{p = 1}^\Gamma \big( 1_{\{\tau_n^+(S_p) \geq S_p + \frac 1 \ell \}} +  1_{\{\tau_n^-(S_p) \leq S_p - \frac 1 \ell \}} \big)  1_{\Omega_n} \big] + o(1),
\end{align}
and we work conditionally on $\mathcal F''_0$ again. We bound $\P \left(\{\tau_n^+(S_p) \geq S_p + \frac 1 \ell \} \cap \Omega_n \big \vert \mathcal F''_0 \right)$ by distinguishing the two cases that the next observation is due to a jump (which has to be the next one after $S_p$ since we are on $\Om_n$) or to an exit of the Brownian motion. In both situations, it is obvious that no exit of $\sigma W_t$ from the interval $[-\epsilon_n c, \epsilon_n c]$  has taken place over $(S_p, S_p + 1/\ell)$. Using Theorem 2.49 in \cite{moeper2010}, this exit time $T$ (starting at $x$) satisfies
\begin{align} \label{moeper}
\E_x[T] = \frac{(\epsilon_n c)^2 - x^2}{\sigma^2} \leq K \epsilon_n^2.
\end{align}
Therefore, 
\begin{align} \label{absch1}
\P \left(\{\tau_n^+(S_p) \geq S_p + \frac 1 \ell \} \cap \Omega_n \big \vert \mathcal F''_0 \right) \leq 2 \sup_{x \in [-\epsilon_n c, \epsilon_n c]} \P_x(T > \frac 1 \ell) \leq {K \epsilon_n^2}{\ell}.  
\end{align}
%
%
%
%
%
%
A similar argument, which will be encountered in the proof of Lemma \ref{step3} in detail (see (\ref{renewal})), can be used to discuss the second indicator in (\ref{abschP}). Thus, for any fixed $\ell$, 
\[
\P \left(\Omega_n(\ell)^c \right) \leq \E[\Gamma] K \epsilon_n^2 \ell + o(1) = o(1).
\]
Note that on $\Omega(\ell) \cap \Omega_n(\ell)$ all jump times and the previous and the following observations are located in the interior of $B_\ell$. Therefore, all intervals with $[\tau_{j-1}^n,\tau_j^n] \cap B_\ell = \emptyset$ are equal in distribution to lengths of intervals of exit times of $\sigma W_t$ from the regular grid. We define
\[
W(\ell)_t = \int_0^t 1_{B_\ell}(s) dW_s,
\] 
and let $({\f''}_t^\ell)_t$ denote the smallest filtration containing $({\f''}_t)_t$ such that $W(\ell)_t$ is ${\f''}_0^\ell$-measurable. We introduce with $Q = Q_\omega$ a regular version of the conditional probability with respect to ${\f''}_0^\ell$. 

We are now in the position to use the stable convergence stated in \cite{fukros2012}, in particular Theorem 2.4 and Proposition 3.1 therein. This is possible, as for any $t$ the convergence $\sup_{i\in\N} |\tau_{i+1}^n\wedge t-\tau_i^n\wedge t| \to 0$ holds in probability, with a similar proof as in their paper. From the fact that $D_n(\ell,t)$ shrinks to $[0,t] \backslash B_\ell$ (using e.g.\ \eqref{absch1} again) we obtain 
\[
1_{\Omega(\ell)\cap \Omega_n(\ell)} \E_{Q_\omega}[ \Psi g(U(\ell)_t^n) ] = 1_{\Omega(\ell)\cap \Omega_n(\ell)} \widetilde \E_{\widetilde Q_\omega}[ \Psi g(\overline U(\ell)_t) ] + o_\P(1),
\]
where $\widetilde Q$ denotes the corresponding probability measure on the product space. Thus, 
\begin{align*}
&\E\big[1_{\Omega_n \cap \Omega(\ell)} \Psi g(U(\ell)_t^n) \prod_{p=1}^k h_p(\alpha(n,p))\big] 
\\ =&\E\big[1_{\Omega(\ell)\cap \Omega_n(\ell)} \Psi g(U(\ell)_t^n) \prod_{p=1}^k h_p(\alpha(n,p))\big] + o_\P(1) \\=& \E\big[1_{\Omega(\ell)\cap \Omega_n(\ell)}\prod_{p=1}^k h_p(\alpha(n,p)) \E_{Q_{\cdot}}[ \Psi g(U(\ell)_t^n) ]\big] + o_\P(1) \\=& \E\big[1_{\Omega(\ell)\cap \Omega_n(\ell)}\prod_{p=1}^k h_p(\alpha(n,p)) \widetilde \E_{\widetilde Q_{\cdot}}[ \Psi g(\overline U(\ell)_t) ]\big] + o_\P(1)  .
\end{align*}
Since $\widetilde \Psi =  \widetilde \E_{\widetilde Q_{\cdot}}[ \Psi g(\overline U(\ell)_t) ]$ is another bounded $\f$-measurable random variable, it is certainly enough to prove 
\begin{align} \label{steplem2}
\E\big[1_{\Omega_n(\ell) \cap \Omega(\ell)} \Xi \prod_{p=1}^k h_p(\alpha(n,p))\big] \to \widetilde \E\big[1_{\Omega(\ell)} \Xi \prod_{p=1}^k h_p(\eta_p)]
\end{align}
for all bounded $\f$-measurable $\Xi$ in order to establish \eqref{step7}. 

%
%
%

On the set $\Omega_n$, the random variables $\alpha(n,p)$ only depend on the jump times $S_p$, $1 \leq p \leq k$, and the Brownian motion over the intervals between those times, as all jumps automatically lead to an observation. By conditioning, it is then sufficient to show (\ref{steplem2}) for
\[
\Xi = \prod_{p=1}^k \varphi_p(\{W_u: u \in (S_{p-1}, S_p]\}) \varrho_p(S_p),
\]
where $\{W_u: u \in (S_{p-1}, S_p]\}$ denotes the Brownian motion between the successive jump times and any $\varphi_p$ and $\varrho_p$ are bounded. Also, a similar argument as the one leading to \eqref{steplem} shows that it is sufficient to work with
\[
\Xi = \prod_{p=1}^k \varphi_p(\{W_u: u \in (S_{p-1}, S_p - 1/m]\}) \varrho_p(S_p),
\]
for any fixed integer $m$ larger than $\ell$. (Note that $S_p - 1/m > S_{p-1}$ always, since we are on $\Omega(\ell)$.) Let us first discuss the asymptotics of
\begin{align*}
& \E\big[h_k(\alpha(n,k)) \varphi_k(\{W_u: u \in (S_{k-1}, S_k - 1/m]\}) \varrho_k(S_k) \vert \f''_{S_{k}-1/m}\big] \\ =& \varphi_k(\{W_u: u \in (S_{k-1}, S_k - 1/m]\}) \varrho_k(S_k) \E\big[h_k(\alpha(n,k))  \vert \f''_{S_{k}-1/m}\big].
\end{align*}
As in (\ref{absch1}), we have
\[
\P(\tau_n^{-}(S_k) \leq S_k - 1/m) \to 0.
\]
On the complement the conditional expectation simplifies, as $\alpha(n,k)$ is independent of $W$ until $S_k - 1/m$ then. We obtain
\[
\E\big[h_k(\alpha(n,k))  \vert \f''_{S_{k}-1/m}\big] = \E[h_k(\alpha(n,k)) \vert \f''_{0}] + o_\P(1).
\]
The idea is to finally use Lemma \ref{step3} below, which proves 
\[
 \E[h_k(\alpha(n,k)) \vert \f''_{0}] \to \widetilde \E[h_k(\eta_k)]. 
\]
One then obtains (\ref{steplem2}) via successive conditioning.  \qed

%
%
%

\begin{lemma} \label{step3}
Let $t > 0$ be arbitrary and let $\tau_n^{-}(t)$ be the last hitting time of the grid $\{\epsilon_n k c~|~ k \in \Z \}$ by $\sigma W$ prior to $t$. Then, as $n \to \infty$ we have the weak convergence
\begin{align*} 
\sigma \epsilon_n^{-1} \left( W_{t} - W_{\tau_n^{-}(t)} \right) \tol c \eta
\end{align*}
where $\eta$ has the density $h$ from (\ref{dens}).

\end{lemma}

\textbf{Proof.} \rm 
Before we begin with the proof, note from self-similarity of Brownian motion that 
$\epsilon_n^{-2} (t- \tau_n^-(t))$ is equal in distribution to $\epsilon_n^{-2}t - \upsilon_n^-(\epsilon_n^{-2}t)$, 
where $\upsilon_n^-(s)$ denotes the last hitting time of the grid $\{k c~|~ k \in \Z \}$ by $\si W_t$. This distribution is called $G_n$ in the following. Also, after a shift to zero, A1.3.0.2 in \cite{borsal2002} shows that the distribution of each hitting time of the grid $\{k c~|~ k \in \Z \}$ by $\si W_t$ equals $S$ with 
\[
\P(S \leq z) = \int_0^{z} \sum_{k=-\infty}^\infty (-1)^k \frac{(2k+1)c}{\sqrt{2 \pi \sigma^2} u^{3/2}} e^{-\frac{((2k+1)c)^2}{2\sigma^2 u}} du = F(z).
\] 
Theorem 1.18 in \cite{mitome2014} then proves 
\begin{align} \label{renewal}
G_n(z) \rightarrow G(z) = \frac{\si^2}{c^2} \int_0^z (1-F(u)) du
\end{align}
pointwise in $z$, where we have used that $c^2/\si^2$ is the expected waiting time until the next hit; recall (\ref{moeper}). 

By conditioning on $t-\tau_n^{-}(t)$ we then obtain 
\begin{align*}
& \P\left( \sigma \epsilon_n^{-1} ( W_{t} - W_{\tau_n^{-}(t)}) \leq x \right) \\ 
=& \int_0^\infty \P\left( \sigma \epsilon_n^{-1} ( W_{t} - W_{\tau_n^{-}(t)}) \leq x \vert t - \tau_n^{-}(t) = \epsilon_n^2 z\right) dG_n(z). 
\end{align*}
Of course, we are not in the setting of stopping times here, as $t - \tau_n^{-}(t) = \epsilon_n^2 z$ states both that $t- \epsilon_n^2 z$ was a hitting time and that no further hitting time has taken place between $t-\epsilon_n^2 z$ and $t$. However, since we work with a regular grid we can at least assume that the place of the last hit was at zero and, using stationarity of the increments of Brownian motion, we can set $t-\epsilon_n^2 z = 0$. Overall, we obtain 
\begin{align*}
& \P\left( \sigma \epsilon_n^{-1} ( W_{t} - W_{\tau_n^{-}(t)}) \leq x \vert t - \tau_n^{-}(t) = \epsilon_n^2 z\right) \\ 
=&  \P\left( \sigma  W_{\epsilon_n^2 z} \leq \epsilon_n x \vert \inf_{0 \leq u \leq \epsilon_n^2 z} \sigma W_u > -\epsilon_n c, \sup_{0 \leq u \leq \epsilon_n^2 z} \sigma W_u < \epsilon_n c \right) \\
=& \P\left( \sigma  W_z \leq x \vert \inf_{0 \leq u \leq z} \sigma W_u > -c, \sup_{0 \leq u \leq z} \sigma W_u < c \right) \\ 
=& \frac{\P\left( \sigma  W_z \leq x, \inf_{0 \leq u \leq z} \sigma W_u > -c, \sup_{0 \leq u \leq z} \sigma W_u < c \right)}{\P\left(\inf_{0 \leq u \leq z} \sigma W_u > - c, \sup_{0 \leq u \leq z} \sigma W_u < c \right)},
\end{align*}
where we have used self-similarity of the Brownian motion again. While the denominator can be written as $\P(S > z)$ with $S$ as before, we use Theorem 33 in \cite{freedman1983} to obtain 
\begin{align*}
& \P\left( \sigma \epsilon_n^{-1} ( W_{t} - W_{\tau_n^{-}(t)}) \leq x \vert t - \tau_n^{-}(t) = \epsilon_n^2 z\right)
\\ =& \frac{\int_{-c}^{x} \frac 1{\sqrt{2 \pi \sigma^2 z}} \sum_{m=-\infty}^\infty (e^{-\frac{(y-4mc)^2}{2 \sigma^2 z}} - e^{-\frac{(y+2c + 4mc)^2}{2 \sigma^2 z}}) dy}{1-F(z)}
\end{align*}
for $x \in [-c,c]$, and it vanishes otherwise. Overall, using dominated convergence, Fubini's theorem and a change of variables,  
\begin{align*}
& \P\left( \sigma \epsilon_n^{-1} ( W_{t} - W_{\tau_n^{-}(t)}) \leq x \right) \\ 
=& \int_0^\infty \P\left( \sigma \epsilon_n^{-1} ( W_{t} - W_{\tau_n^{-}(t)}) \leq x \vert t - \tau_n^{-}(t) = \epsilon_n^2 z\right) dG_n(z) \\
\to& \int_0^\infty \frac{\int_{-c}^{x} \frac 1{\sqrt{2 \pi \sigma^2 z}} \sum_{m=-\infty}^\infty (e^{-\frac{(y-4mc)^2}{2 \sigma^2 z}} - e^{-\frac{(y+2c + 4mc)^2}{2 \sigma^2 z}}) dy}{1-F(z)} dG(z) 1_{[-c,c]}(x) \\
=& \int_{-c}^{x} \int_0^\infty \frac 1{\sqrt{2 \pi c^2 u}} \sum_{m=-\infty}^\infty (e^{-\frac{(\frac yc-4m)^2}{2u}} - e^{-\frac{(\frac yc+2 + 4m)^2}{2 u}}) du dy 1_{[-c,c]}(x) \\ =& \int_{-\infty}^{x} \frac 1c h\Big(\frac yc \Big) dy. 
\end{align*}
This proves the result. Note finally that $h$ is easily seen to be a density, as 
\begin{align*}
\int_{-1}^{1} h(y) dy = \int_0^\infty \P\left(\inf_{0 \leq u \leq z} W_u > -1, \sup_{0 \leq u \leq z} W_u < 1 \right) dz = 1
\end{align*}
from $\int_0^\infty \P(X > z)dz  = \E[X]$ for any positive random variable and (\ref{moeper}). \qed \\

Let us now finish the proof of Theorem \ref{thm1} when $X_t = \si W_t + J_t$. From (\ref{moeper}) we obtain 
\[
\epsilon_n^{-1} \left([X,X]_t  - [X,X]_{\tau_n^{-}(t)} \right) = O_\P(\epsilon_n) = o_\P(1),
\] 
because $t - \tau_n^{-}(t)$ is smaller in distribution than $T$. As before, we then have the key decomposition
\begin{align}
& \epsilon_n^{-1} \left( RV(X,X)_t^n - [X,X]_{\tau_n^{-}(t)} \right) \nonumber \\ \label{decomp2}
=&  U_t^n + \epsilon_n^{-1} \left( RV(J,J)_t^n - [J,J]_{\tau_n^{-}(t)} \right) + 2 \epsilon_n^{-1} RV(C,J)_t^n,
\end{align}
and for simplicitly we can assume to be on $\Omega_n$ at the cost of additional smaller order terms only. On this set every jump causes an observation and the second term in \eqref{decomp2} vanishes identically. The third term simplifies as well, as we end up with 
\begin{align*}
\widetilde Z_t^n(X) &= \epsilon_n^{-1} \left( RV(X,X)_t^n - [X,X]_t \right) \\ &= U_t^n + 2 \epsilon_n^{-1} \sum_{S_p \leq t} \Delta J_{S_p} \sigma \left( W_{S_p} - W_{\tau_n^{-}(S_p)} \right) + o_\P(1).
\end{align*}
The result then follows from Lemma \ref{lemstab}, as $\overline U_t$ equals $\widetilde U_t$ for a constant $\si$. 

Finally, let us explain why it is sufficient to discuss the simple case 
\begin{align*}
\hat X_t = \sigma W_t + J_t
\end{align*}
only. To this end, let 
\begin{align*}
X_t = \int_0^t b_s ds + \int_0^t \sigma_s dW_s + J_t
\end{align*}
as above. Following Theorem IV.4.32 in \cite{jacshi2003}, there exists an equivalent probability measure $\Q$ such that the continuous martingale part remains the same and the jump process is still of finite activity under $\Q$, but the drift vanishes. Note further that the $\mathcal F$-conditional distribution of the limiting process $\widetilde Z_t = \widetilde{U}_t + \widetilde{V}_t$ only depends on the quadratic variation process $[X,X]$ up to time $t$ which remains unchanged under a change of measure using Theorem III.3.13 in \cite{jacshi2003}. Thus, if we can show $\mathcal F$-stable convergence under $\Q$, we directly get
\begin{align*}
\E_\P[Y g(\widetilde Z_t^n(X))] &= \E_\Q\big[Y \frac{d\P}{d\Q} g\big(\widetilde Z_t^n(X)\big)\big] \\ &\rightarrow \E_\Q\big[Y \frac{d\P}{d\Q} g(\widetilde Z_t)\big] = \E_\P\big[Y g(\widetilde Z_t)\big]
\end{align*}
for any bounded, $\F$-measurable $Y$ and any bounded and continuous function $g$, as $\E_\Q[d\P/d\Q|\F_t]$ is locally bounded due to Proposition III.3.5 in \cite{jacshi2003} and thus can be assumed to be bounded as well.  

Suppose therefore that 
\begin{align*}
X_t = \int_0^t \sigma_s dW_s + J_t.
\end{align*}
Using a time change, we can write $X = \hat X_{[X^c,X^c]}$ with $\hat X_t = \hat C_t + \hat J_t,$ where $\hat C$ is a standard Brownian motion and $\hat J$ is a finite activity jump process with jump times $\hat S_p = [X^c,X^c]^{-1}(S_p)$ and the same jump sizes. We then know from the preceding results that 
\[
\epsilon_n^{-1} \left( RV(\hat X,\hat X)_t^n - [\hat X,\hat X]_t \right) \tols \sqrt{\frac 23} c \int_0^t dW'_s + 2 \sum_{\hat S_p \leq t} \Delta \hat X_{\hat S_p} c \eta_p.
\] 
The general claim then follows from 
\[
RV(\hat X,\hat X)_{[X^c, X^c]_t}^n - [\hat X,\hat X]_{[X^c, X^c]_t} = RV(X,X)_{t}^n - [X,X]_{t}
\]
using the fact that the exit times commute with the specific choice of the time change, as well as from 
\[
\int_0^{[X^c,X^c]_t} dW'_s = \int_0^{t} dW'_{[X^c,X^c]_s} = \int_0^t \sigma_s dW'_s 
\]
and 
\[
\sum_{\hat S_p \leq [X^c,X^c]_t} \Delta \hat X_{\hat S_p} c \eta_p = \sum_{S_p \leq t} \Delta X_{S_p} c \eta_p.
\]
\qed

\bibliographystyle{chicago}
\bibliography{biblio}
\end{document}